\documentclass[letterpaper, 10 pt, conference]{ieeeconf}  

\IEEEoverridecommandlockouts                              

\usepackage{graphics} 
\usepackage{epsfig} 
\usepackage{mathptmx} 
\usepackage{times} 
\usepackage{amsmath} 
\usepackage{amssymb}  
\usepackage[margin=3cm]{geometry}
\usepackage{algorithm2e}
\usepackage{comment}

\title{\LARGE \bf
Proximal Gradient-based Low Rank Tensor Decomposition for State Dependent Riccati Equation*}

\author{Jiahua Jiang$^{1}$ and Carmeliza Navasca$^{2}$
\thanks{*This work has been submitted to the IEEE for possible publication. Copyright may be transferred without notice, after which this version may no longer be accessible.}
\thanks{$^{1}$Jiahua Jiang is with the Institute of Mathematical Science,
        ShanghaiTech University, Shanghai,  China 
        {\tt\small jiangjh@shanghaitech.edu.cn}}%
\thanks{$^{2}$Carmeliza Navasca is with the Department of Mathematics, University of Alabama at Birmingham,
        Birmingham, AL 35294, USA
        {\tt\small cnavasca@uab.edu}}%
}

\begin{document}

\maketitle
\thispagestyle{empty}
\pagestyle{empty}

\begin{abstract}

We address the optimal control problems arising from partial differential equations with large discrete dimensional control systems. To obtain reduced order models, we find basis elements from the canonical polyadic (CP) decomposition. Tensor datasets are from snapshots of the large models. Our method to reduce the control system is to use dimensionality reduction approaches through sparse optimization and flexible hybrid methods is to obtain low rank CP tensor basis elements. The reduced optimal control problem leads to reduced state-dependent Riccati Equations which can be solved efficiently.

\end{abstract}

\section{INTRODUCTION}

Control strategies play an important role in many engineering applications and most recently, in AI technology. Some of these applications are in autonomous vehicles, semiconductor device fabrication, fluid flow and deep learning. An nonlinear optimal  control problem formulation is a standard way to find nonlinear control feedback strategies. The solutions to the optimization satisfy the well known Hamilton-Jacobi-Bellman (HJB) Equations. It is well known that the computational cost is prohibitively expensive when standard methods are implemented. In this work, we will be constructing nonlinear control feedback through solving State-Dependent Riccati Equation (SDRE). SDRE was introduced by \cite{Pearson} in 1962 but there have been many advancements, e.g. see \cite{Crandall-Lions}.

Our proposed method is to use tensor computation for problems in control systems. In recent years, low rank tensor decompositions have been a vital tool in data science and engineering disciplines. It is currently playing a major role in deep neural network applications and architectures. There are already many contributions of tensor methods in control systems. Here are some recent developments of tensor methods in control problems.
In \cite{Heidel-Khoromskaia-Khoromskij-Schulz}\cite{Schmitt-Khoromskij-Khoromskaia}, discretized optimal control problems with 2D and 3D elliptic operators and coefficients are solved using tensor methods. More specifically, they have shown low rank preconditioned conjugate gradient iterative technique is effective in solving the Karush-Kuhn-Tucker conditions of the tensor-structured lagrangian formulation. In \cite{Kirsten-Saluzzi}, through the dynamic programming formulation on a tree structure of 
the semi-discrete finite optimal control problem with  semilinear PDEs as contraints is derived. It is solved through a reduced order model, derived by applying the Higher-Order Orthogonal Iteration (HOOI), aka higher-order SVD (HO-SVD), on the tensor structured snapshots, where the $i$th snaphot is the $i$th slice of the PDE solution at time $t_i$.  
In \cite{Borggaard-Zietsman}, the quadratic-quadratic regulator problem is solved by expanding the cost and dynamics functions in their Taylor series expansions where series of tensor equations are solved efficiently by utilizing special properties of Kronecker products.
In \cite{Dolgov-Kalise-Kunisch}, the value function is represented by separable basis which leads to a system of Galerkin equations that is solved through low-rank tensor train decomposition.

The main contributions of this paper are as follows.
\begin{itemize}
    
    \item The proximal gradient method is applied to nonlinear control datasets to find sharp tensor rank while using flexible hybrid methods for regularization.
    \item The tensor rank $R$ dictates the reduced order dimension of the state system.
    \item Our method leads to algorithms that find reduced basis elements for the control systems and construct control feedback which stabilizes the state dynamics.
   
\end{itemize}
The advantages of our methods are as follows.
\begin{itemize}
    \item  Tensor data structures are  easily constructed as frontal slices $\mathcal{T}(:,:,k)$ of the solution profiles for $k$ control samples.
     \item 
    The approximated control feedback exhibits fast stabilization of the control system.
    \item The on-line computational complexities of the proposed method are much smaller than the computational cost $\mathcal{O}(N_tN_x^3)$ of the state dependent Riccati equation. 
    \item The computed objective cost function associated the control feedback from the reduced models are a small fraction of the cost of the full model $J_{full}$.
\end{itemize}

\section{Problem Formulation}
Given a time varying ordinary differential equation,
\begin{eqnarray}\label{mainode}
\dot{x}(t)&=&A(x(t))x(t) + B(x(t))u(t), \\
x(0)&=&x_0, \nonumber
\end{eqnarray}
the goal is to drive the system $x(t) \longrightarrow 0$ as $t \longrightarrow \infty$ where $u(t):[0,\infty) \rightarrow \mathbb{R}^m$ is unknown control function. Let 
$x(t):[0,\infty) \rightarrow \mathbb{R}^d$ be the state of the system, $A(t) \in \mathbb{R}^{d \times d}$ and $B(t) \in \mathbb{R}^{d \times m}$ are state dependent matrices. To find the optimal unknown control $u$, an optimization is formulated,
\begin{equation}\label{mainopt}
\begin{aligned}
& \underset{u \in \mathcal{U}}{\text{minimize}}
& &  \int_0^{\infty} \Vert x(t) \Vert^2_Q + \Vert u(t) \Vert^2_R~dt \\
& \text{subject to}
& & \dot{x}(t)=A(x(t))x(t) + B(x(t))u(t), ~ t \in (0,\infty), \\
&&& x(0) = x_0
\end{aligned}
\end{equation}
where $\Vert x(t) \Vert^2_Q = x(t)^TQx(t)$ and 
$\Vert u(t) \Vert^2_R  = u^TRu$. Often, $Q$ and $R$ are assumed  to be positive definite; i.e. $Q_{ij}=Q_{ji}$, $R_{ij}=R_{ji}$ and the eigenvalues of $Q$ and $R$ are positive.


\subsection{Optimality Conditions}
Let's consider a general problem, the infinite-time-horizon optimal control problem,
    \begin{equation}
    \begin{aligned}\label{eqn:optContrCost}
    \min_{u} \int_0^\infty \ell(x(t),u(t)) dt
    \end{aligned}
    \end{equation}
subject to the controlled dynamics
    \begin{equation}
    \begin{aligned}
        \dot{x}(t) &=f(x(t),u(t)) \label{eqn:dynamics}\\
        x(0) &= x_0.
    \end{aligned}
    \end{equation}
    where
    \begin{itemize}
        \item $x: [s,\infty) \to \mathbb{R}^n$ is the state as a function of time
        \item $u: [s,\infty) \to \mathbb{R}^m$ is the control or trajectory as a function of time
        \item $\ell: \mathbb{R}^n \times \mathbb{R}^m \to \mathbb{R}$ is the Lagrangian or running cost and is assumed to be smooth and strictly convex in $u$
        \item $f: \mathbb{R}^n \times \mathbb{R}^m \to \mathbb{R}^n$ is the dynamics function is and is assumed to be smooth.
    \end{itemize}

Define a value function,
\begin{eqnarray}
    J(u(x))= \int_0^{\infty} l(x,u)~dt
\end{eqnarray}
that satisfies
\begin{equation}\label{hjb1}
\frac{\partial J}{\partial x}(x)f(x,u^*) + l(x,u^*)=0
\end{equation}
where we denote $u^*$ as the optimal control solution.
The derivation of (\ref{hjb1}) is
through the notion of viscosity solution; see e.g. \cite{Crandall-Lions}, \cite{Navasca-Krener1}, \cite{Navasca-Krener2}. Moreover, the optimal control $u^*$ is found through the optimization:
\begin{equation}\label{hjb2a}
   \min_u \left [ \frac{\partial J}{\partial x}(x)f(x,u) + l(x,u) \right ]
\end{equation}
which is equivalent to
\begin{eqnarray} \label{hjb2b}
    \frac{\partial J}{\partial x}(x)\frac{\partial f}{\partial u}(x,u^*) + \frac{\partial l}{\partial u}(x,u^*)=0
\end{eqnarray}
Equations (\ref{hjb1}-\ref{hjb2a}) are the optimality conditions in which the unknown solutions $u^*$ (optimal control) and $J(u^*(x))$ (optimal cost) satisfy. These are called the Hamilton-Jacobi-Bellman PDEs.

\subsection{State Dependent Riccati Equation}
Recall from (\ref{mainopt}) that we have 
\begin{equation*}
    \ell(x(t),u(t))=\Vert x(t) \Vert^2_Q + \Vert u(t) \Vert^2_R
\end{equation*}
\begin{center}
and 
\end{center}
\begin{equation*}
    f(x(t),u(t))=A(x(t))x(t) + B(x(t))u(t).
\end{equation*}
 We can assume the solution is of the form $J(u^*(x))=(x(t))^T\Pi(x(t)) x(t)$ where $\Pi$ is an $d \times d$ matrix to be approximated. Similarly, we can assume $u^*(x)=K(x(t))x$ where $K$ is an $m \times d$ unknown matrix. It has been shown in \cite{Anderson-Moore} that these assumptions are the correct form. \\

The optimality condition (\ref{hjb2a}) for control $u$ is
$$2x^T(\Pi(x(t)) B + K(x(t))^TR)=0$$
which implies $K=-R^{-1}B^T \Pi(x(t))$ and $u^*(x)=Kx=-R^{-1}B^T \Pi(x(t))x.$ Recall that the following assumptions:  $Q=Q^T,~R=R^T$. By plugging in $u^*(x)=-R^{-1}B^T \Pi(x(t)) x$ and $\frac{\partial J}{\partial x}(x)= 2 (x(t))^T\Pi(x(t))$ into the optimality condition (\ref{hjb1}), we obtain
the optimality condition 
\begin{equation}
    \begin{aligned}
        x^T(\boxed{A^\top \Pi + \Pi A + Q - \Pi B R^{-1} B^\top \Pi})x = 0
    \end{aligned}
\end{equation} 
By equating the boxed terms to zero, we arrive at the state dependent Riccati equation (SDRE).
\begin{eqnarray}
    &&A(x(t))^\top \Pi(x(t)) + \Pi(x(t)) A(x(t))\\ \nonumber
    &&+ Q - \Pi(x(t)) B R^{-1} B^\top \Pi(x(t))=0
\end{eqnarray} \label{sdre}

Finally, the approximated matrix $K$ must satisfy must satisfy the stability condition, that is, the eigenvalues of $A(x(t))+BK(x(t))$ must be negative or have negative real part for complex eigenvalues, $Re(\lambda (A+BK) )< 0$.
This stability condition implies that $x(t) \rightarrow 0$ as $t \rightarrow \infty$.

\section{Low Rank Tensor Decomposition}
We consider a given tensor data $\mathcal{T} \in \mathbb{R}^{I \times J \times K}$  to find the \emph{closest} low rank tensor using the canonical polyadic (CP ) decomposition:
\begin{equation}\label{cpd2}
\mathcal{T} \approx \sum_{r=1}^{R} \alpha_r\mathbf{x}_r\circ\mathbf{y}_r\circ\mathbf{z}_r
\end{equation}
where $\alpha_r\in\mathbb{R}$ is a rescaling coefficient of rank-one third-order tensor $\mathbf{x}_r\circ\mathbf{y}_r\circ\mathbf{z}_r$, an outer product of three vectors, for $r=1,\cdots,R$. Let ${\alpha}=(\alpha_1,\cdots,\alpha_R)\in\mathbb{R}^R$
and denote $[{\alpha};\mathbf{X},\mathbf{Y},\mathbf{Z}]_R = \sum_{r=1}^{R} \alpha_r\mathbf{x}_r\circ\mathbf{y}_r\circ\mathbf{z}_r$
in (\ref{cpd2}) with $\mathbf{X}=(\mathbf{x}_1,\cdots,\mathbf{x}_R)\in\mathbb{R}^{I\times R},\mathbf{Y}=(\mathbf{y}_1,\cdots,\mathbf{y}_R)\in\mathbb{R}^{J\times R}$ and $\mathbf{Z}=(\mathbf{z}_1,\cdots,\mathbf{z}_R)\in\mathbb{R}^{K\times R}$
are called the factor matrices of a third-order tensor $\mathcal{T}$.

In most iterative techniques for decompositions, tensor \emph{matricizations} are required to transform the tensor equations into matrix equations. Using the Khatri-Rao product, a column-wise kronecker product, denoted as $\odot$, the CP model (\ref{cpd2}) can be written in three equivalent matrix equations: $\mathbf{T}_{(1)}
  = \mathbf{XD}(\mathbf{Z\odot Y})^T$, 
  $\mathbf{T}_{(2)}
  =\mathbf{YD}(\mathbf{Z\odot X})^T$  and 
  $\mathbf{T}_{(3)}
  =\mathbf{ZD}(\mathbf{Y\odot X})^T$ where the matrix $\mathbf{D}$ is diagonal with elements of $\alpha$.  To reconstruct  $\mathcal{T}$ in a CP decomposition fashion with a known tensor rank $R$ and $\mathbf{D}=\mathbf{I}$, the linear least-squares subproblems are implemented to solve iteratively for $\mathbf{X}$, $\mathbf{Y}$ and $\mathbf{Z}$. 
see \cite{Navasca-DeLathauwer-Kindermann}, \cite{Li-Kindermann-Navasca}, \cite{Glenn-Li-Navasca}, \cite{Dulal-Karim-Navasca} 
for Alternating Least-Squares (ALS) method. See Algorithm \ref{alg:als}.
\subsection{Tensor Rank via Sparse Optimization}
To approximate the tensor $R$ and reconstruct $\mathcal{T}$ into a CP format, we formulate a sparse optimization problem \cite{Wang-Navasca} for recovering CP decomposition from tensor $\mathcal{T} \in \mathbb{R}^{I \times J \times K}$: 
\begin{subequations}
\begin{align}
\min_{A,B,C,\sigma} \left \Vert \mathcal{T} - \mathcal{S} \right \Vert_F + \lambda \Vert \alpha \Vert_{\ell_1} \label{prob}
  \end{align}
\end{subequations}
where $\mathcal{}$ $\lambda$ is a constant regularization parameter and  $\mathcal{S}= \sum_r \alpha_r    
\mathbf{x}_r\circ\mathbf{y}_r\circ\mathbf{z}_r$. The ALS in  Algorithm \ref{alg:als} along with 
\begin{equation}\label{alphahOPT}
\alpha \longleftarrow \min \limits_{\alpha} \frac{1}{2} \|\mathbf{t}-{\alpha^T \mathbf{Q}}\|_F^2 + \lambda \Vert \alpha \Vert_1
\end{equation}
where $\mathbf{t}$ and $\mathbf{Q}$ are vectorization of $\mathcal{T}$ and $\mathbf{x}_r\circ\mathbf{y}_r\circ\mathbf{z}_r$, respectively, are the updating schemes for $\mathbf{X}, \mathbf{Y}, \mathbf{Z}$ and $\alpha$.
Note that the gradient of $\frac{1}{2} \Vert \mathcal{T} - \mathbf{x}_r\circ\mathbf{y}_r\circ\mathbf{z}_r \Vert_F^2 $ is 
 $\nabla_{\alpha} f(\mathbf{X,Y,Z,}{\alpha})=({\alpha}^T\mathbf{Q}-\mathbf{t})\mathbf{Q}^T$ which leads to the least-squares formulation in (\ref{alphahOPT}).

 The equation (\ref{alphahOPT}) is associated with Iterative Soft Thresholding Algorithm (ISTA) \cite{Wang-Navasca} \cite{Jiang-Navasca-Sanogo} whose derivation is based on the Majorization-Minimization method.    

Choosing the regularization parameter $\lambda$  plays a crucial role in solving (\ref{alphahOPT}). The proximal operators formulation and the MM approach are used to solve $\alpha$ iteratively via 
\begin{eqnarray*}
\alpha \leftarrow \min\limits_{{\alpha}} \{\langle{\alpha}-{\alpha}^n,
\nabla_{{\alpha}}f\rangle
+\frac{s\eta_n}{2}\|{\alpha}-{\alpha}^n\|^2+\lambda\|{\alpha}\|_1\}.
\end{eqnarray*}
In \cite{Wang-Navasca}, theoretical bounds for the regularization parameters are provided to help with parameter tuning. However, the bounds can be impractical at times. To address this problem, we embedded the flexible hybrid method  into the CP decomposition framework; see e.g. \cite{Jiang-Navasca-Sanogo}. The idea is to replace $\lambda \Vert  \alpha \Vert_{\ell_1} \approx \lambda \Vert L_k \alpha \Vert_{\ell_2}$ sequentially for some $L_k$. This leads to a sequence of optimization problems:
$$
\min_{\hat{\alpha}} \Vert (L^{-1}_k \hat{\alpha})^T \mathbf{Q}-\mathbf{t} \Vert^2_{\ell_2} + \lambda \Vert  \hat{\alpha}\Vert_{\ell_2}$$ 
where $\hat{\alpha}=L_k \alpha$. Moreover, the Golub-Kahan bidiagonalization is applied to further reduced to
$$
 \min_y \Vert M_k y - \beta_1 e_1 \Vert_{\ell_2} + \lambda \Vert y \Vert_{\ell_2}
$$
satisfying $\mathbf{Q}^TZ_k=U_{k+1}M_k$ and $\mathbf{Q}U_{k+1}=V_{k+1}T_{k+1}$ where $M_k$ and $T_k$ are upper Hessenberg matrices, $U_{k+1}$ and $V_{k+1}$ have orthogonal columns and $Z_k=L^{-1}_k v_k$ such that $\alpha_k=Z_k y_k$.

\section{Numerical Approach}

\subsection{Uncontrolled State Dynamics}

The Allen Cahn model \cite{ZhangDu} is the following,
\begin{eqnarray}\label{1dAC}
    \frac{dv}{dt}(x)= \nu v_{xx} - \frac{1}{2\xi^2} v(x)(v(x)^2 -1)
\end{eqnarray}
where the time and state domains are $t \in [0,1]$ and $x \in [0,2]$, respectively, with boundary condition, $\frac{\partial v}{\partial n}=0$. The boundary condition is handled numerically as $v(1)=v(2)$ and $v(nx)=v(nx-1)$. The initial condition is generated by using cubic spline interpolant to fit with a dataset in \cite{Morosanu}. The parameters are $\nu=1$ and $\xi=0.02$. Figure \ref{figure0} depicts the numerical solution of the Allen-Cahn equation (\ref{1dAC}).

\subsection{Construction of Tensor Data}
The raw tensor data is of the size $101 \times  550 \times 50$ with modal units: spacial discretization $\times$ time discretization $\times$ $\beta$ values. The $\beta$ values were selected uniformly at random from the interval $[-100,100]$. The frontal slice in the $k$th mode is an $101 \times 550 $ matrix, $\mathcal{T}(:,:,k)$ from the solution profile of 
\begin{eqnarray}\label{1dACb}
    \frac{dv}{dt}(x)= \nu v_{xx} - \frac{1}{2\xi^2} v(x)(v(x)^2 -1) + u(x)
\end{eqnarray}
where $u(x)=\beta_i v(x)$ with 
$\beta_i \in[-100,100]$.
\subsection{Tensor Decomposition}
Here we describe the algorithms, Proximal Gradient for Sparsity (PGS) and ALS, for CP decompositions in Algorithm \ref{alg:als} and Algorithm \ref{alg:pgs}. We combined the two algorthms, to form PGS+ALSk where tensor rank $R$ is calculated from PGS while the factor matrices $\mathbf{X,Y,Z}$ are approximated using ALS. Note that PGS+ALSk is when the input tensor rank is $R+(k-1)$ in ALS with $k-1$ added basis elements.

\begin{algorithm}
\caption{Alternating Least-Squares (ALS)}\label{alg:als}
\KwData{third-order tensor $\mathcal{T}$  }
\KwData{tensor rank $R$  }
\KwData{initial factors $\mathbf{X}_0$, $\mathbf{Y}_0$,$\mathbf{Z}_0$ }
\KwResult{tensor rank $R$ and factors: $\mathbf{X}$, $\mathbf{Y}$,$\mathbf{Z}$}
stopping criteria=$\epsilon$\;

 \While{$err>\epsilon$}{
  $\mathbf{X}_i \longleftarrow \min_\mathbf{X} \Vert  \mathbf{T}_{(1)}
  - \mathbf{X}(\mathbf{Z_i\odot Y_i})^T \Vert$ \;
  $ \mathbf{Y}_i \longleftarrow \min_\mathbf{Y} \Vert \mathbf{T}_{(2)} - \mathbf{Y}_i(\mathbf{Z_i\odot X_{i+1}})^T \Vert$ \;
  $\mathbf{Z}_i \longleftarrow \min_\mathbf{Z} \Vert  \mathbf{T}_{(3)} - \mathbf{Z}(\mathbf{Y_{i+1}\odot X_{i+1}})^T \Vert $\;
  $err=\Vert \mathcal{T}- \sum_{r=1}^{R} \mathbf{X_i(:,r)}\circ\mathbf{Y_i(:,r)}\circ\mathbf{Z_i(:,r)} \Vert $
  }
 \end{algorithm}

\begin{algorithm}
\caption{Proximal Gradient- Sparse (PGS)}\label{alg:pgs}
\KwData{third-order $I\times K \times K$ tensor $\mathcal{T}$  }
\KwData{upperbound tensor rank $\hat{R}=min\{IJ, IK, JK\}$  }
\KwData{initial factors $\mathbf{X}_0$, $\mathbf{Y}_0$,$\mathbf{Z}_0$ }
\KwData{coefficient $\mathbf{\alpha_0}$}
\KwResult{tensor rank $R$}
\KwResult{ $\mathbf{X}$, $\mathbf{Y}$,$\mathbf{Z}$}
\KwResult{$\mathbf{\alpha}$}
stopping criteria=$\epsilon$\;

 \While{$err>\epsilon$}{
 $\mathbf{D}=diag(\mathbf{\alpha})$ \;
 $ \mathbf{X} \longleftarrow \min_\mathbf{X} \Vert  \mathbf{T}_{(1)}
  - \mathbf{XD}(\mathbf{Z\odot Y})^T \Vert$ \;
 $\mathbf{Y} \longleftarrow \min_\mathbf{Y} \Vert \mathbf{T}_{(2)} - \mathbf{YD}(\mathbf{Z\odot X})^T \Vert$\;
   $\mathbf{Z} \longleftarrow \min_\mathbf{Z} \Vert  \mathbf{T}_{(3)} - \mathbf{ZD}(\mathbf{Y\odot X})^T \Vert$ \;
  $\alpha \longleftarrow \min \limits_{\alpha} \frac{1}{2} \|\mathbf{t}-{\alpha^T \mathbf{Q}}\|_F^2 + \lambda \Vert \alpha \Vert_1$ \;
  $err=\Vert \mathcal{T}- \sum_{r=1}^{\hat{R}} \mathbf{\alpha_r}\mathbf{x_r}\circ\mathbf{y_r}\circ\mathbf{z_r} \Vert $
  }
$R=\Vert \mathbf{\alpha} \Vert_0$ 
 \end{algorithm}

\subsection{Reduced Order Model}
The standard discretization of Allen-Cahn equation (\ref{1dACb}) leads to the full order model,
\begin{eqnarray} \label{linear_system_1d_a}
\dot{v}=\underbrace{\left (\nu \Delta_d + \frac{1}{2\xi^2} I_d - \frac{1}{2\xi^2} D_d  \right )}_Av + Bv
\end{eqnarray}
where $x_i \in [0,2]$ is the $i$th uniformly discretized point for $i=1,\ldots,d$ and $v_i=v(x_i)$ with $d=101$. The matrices $I_d$ and $D_d$ are the $d$-dimensional identity matrix and the $d$-dimensional diagonal matrix with $v^2_i$ on its diagonal entries for $i=1,\ldots,d$, i.e.
\begin{eqnarray*}
    D_d=\left [ \begin{array}{cccc}
    v^2_1 & 0   & \ldots & 0 \\
    0   & v^2_2 & \ldots & 0 \\
    0   & 0   & \ddots & 0 \\
    0   & 0   & \ldots & v^2_d 
    \end{array}
    \right].
\end{eqnarray*}
 Also, $B=-\beta I_d$ to generate the tensor data $\mathcal{T}$.

An orthogonal matrix $P \in \mathbb{R}^{101 \times nxred}$ such that $v=Pv_{red}$ is obtained from a low rank reconstruction of $\mathcal{T}$ through CP decomposition. The columns of $P$ are the left singular values of the product of the factor matrices of $\mathcal{T}$, $\mathbf{X}\mathbf{Y}^T$. These are steps in the implementation of POD and reduced basis methods, e.g. \cite{Chen-Gottlieb}\cite{Jiang-Navasca-Sanogo}.The reduced order model is derived as follows:
 \begin{eqnarray*}
    \dot{v}_{red}&=&A_{red}v_{red} + B_{red}u
    \end{eqnarray*}
    where $A_{red}=P^TAP $ and $B_{red}=P^TB$.

    \begin{algorithm}
\caption{SDRE for reduced order (SDRERO)}\label{alg:sdrero}
\KwData{$\dot{x}=A_t x+B_t u$,  }
\KwData{$l(x,u)=x^TQ_tx + u^TR_tu$}
\KwResult{control $u^*(x)$, value $J(u^*(x))$}
$t=linspace(0,1,Nt)$\;
$\Delta = \frac{1}{Nt-1}$\;
 \For{i=1:Nt-1}{
  $\Pi_{i}  \longleftarrow \Pi_{i-1} A_{i-1} + A_{i-1}^T + Q_{i-1} + \Pi_{i-1} B_{i-1} R^{-1}_{i-1}B_{i-1}=0$\;
  $K_{i} \longleftarrow -R_{i-1}B^T_{i-1}\Pi_{i-1}$\;
  $u_{i} \longleftarrow K_{i}x_i$\;
  $x_{i+1} \longleftarrow (1 + \Delta t (A_i+B_iK_i) ) x_{i}$\;
  }
 $u^*_{Nt}=K_{Nt-1}x_{Nt}$\;
 $J=x^T_{Nt}\Pi_{Nt-1}x_{Nt}$
 \end{algorithm}
The reduced optimal control problem is  solve through SDRE.
Although this approach, which can seen in Algorithm $1$, is prohibitively expensive for Riccati equations since the computational complexity is $\mathcal{O}(d^3)$ per time step (i.e. Bartel-Stewart \cite{Bartel-Stewart}). However, this is proved to be effective for smaller order dimension such as, $k<<d$. Thus, this method is highly effective for the reduced order models. Similar approach can be seen in \cite{Banks}. In addition, Algorithm $1$ can be accelerated by the Hessenberg-Schur decomposition \cite{GolubNashVanLoan} and other iterative methods like biconjugate gradient descent \cite{BrazellLiNavascaTamon}.

\section{Numerical Results}

\begin{figure}[thpb]
      \centering
      \framebox{\parbox{2.75in}{ 
       \begin{center}
      \includegraphics[scale=.35]{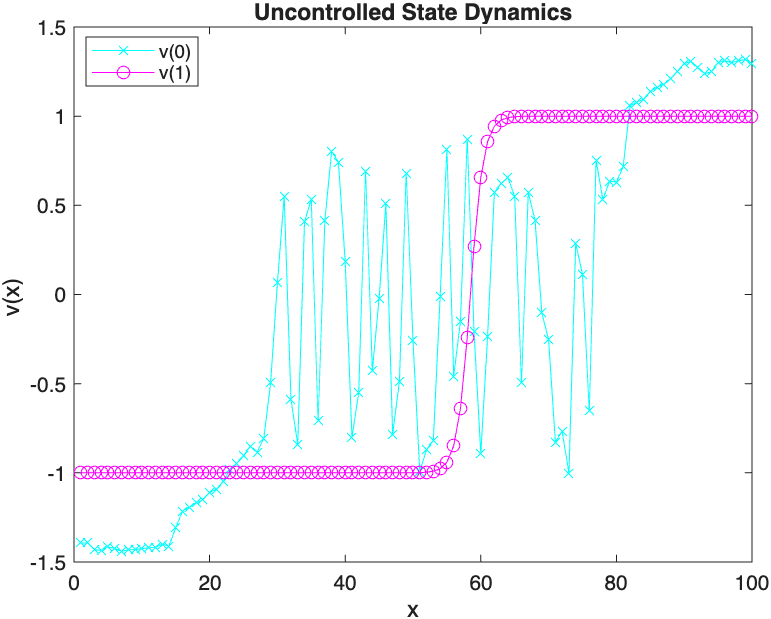}
      \end{center}
      \caption{Uncontrolled state dynamics and initial condition.}
      \label{figure0}
      }}
   \end{figure}

In our experiments, we have $A \in \mathbb{R}^{101 \times 101}$ and $B \in \mathbb{R}^{101 \times 2}$ in the full order model (\ref{linear_system_1d_a}). It follows that the state dynamics $v \in \mathbb{R}^{101}$ and the control variable $u \in \mathbb{R}^2$.  We implement the algorithms PGS and PGS+ALSk where $k=1,2$ (Algorithms \ref{alg:als}-\ref{alg:pgs}) to find $A_{red}$ and $B_{red}$. Here we found $R=2$. The initial condition is based on the interpolation of data points in \cite{Morosanu}. Subsequently, we set the initial condition as $v_{red}(0)=Pv(0)$. See Figure \ref{figure0} for the initial condition graph.

\begin{figure}[thpb]
      \centering
      \framebox{\parbox{2.75in}{ 
       \begin{center}
      \includegraphics[scale=.55]{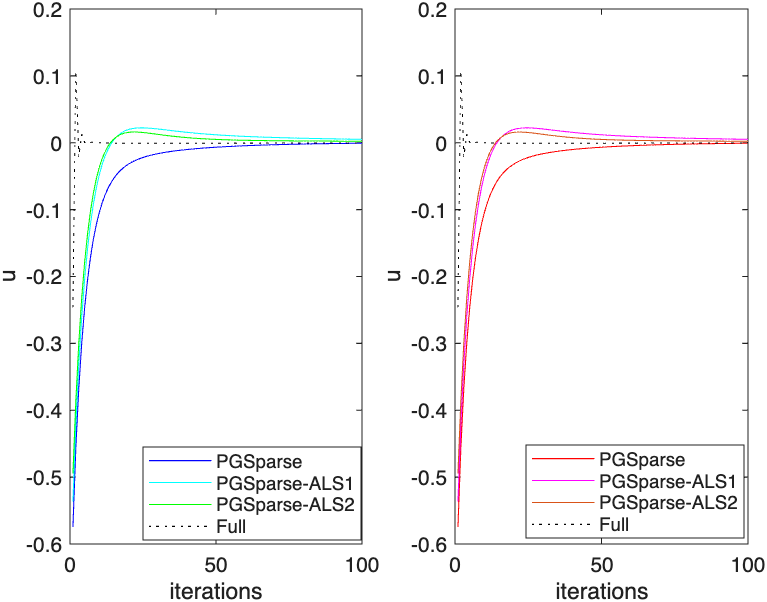}
      \end{center}
      \caption{Feedback control.}
      \label{figure1a}
      }}
   \end{figure}

   \begin{figure}[thpb]
      \centering
      \framebox{\parbox{2.85in}{ 
       \begin{center}
      \includegraphics[scale=.57]{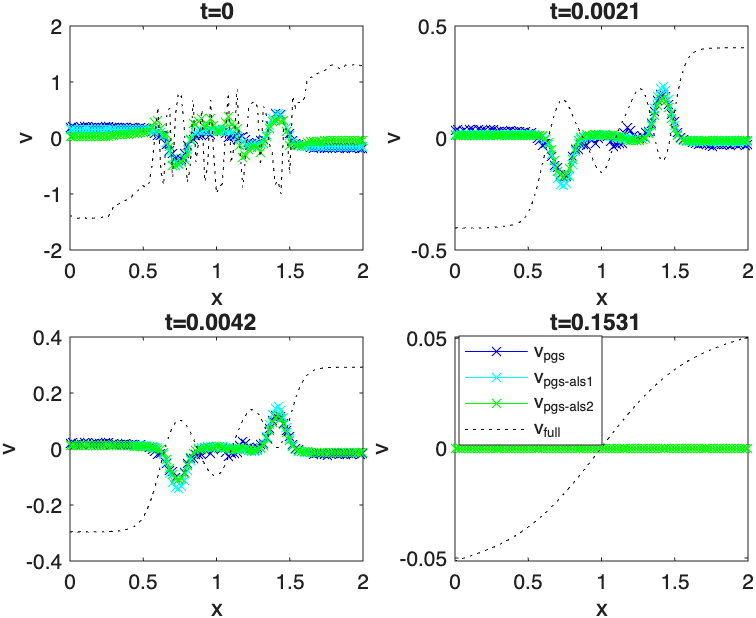}
        \end{center}
      \caption{State with Feedback Control}\label{figure1}
      }}
      \end{figure}

    \begin{figure}[thpb]
      \centering
      \framebox{\parbox{2.85in}{
       \begin{center}
      \includegraphics[scale=.5]{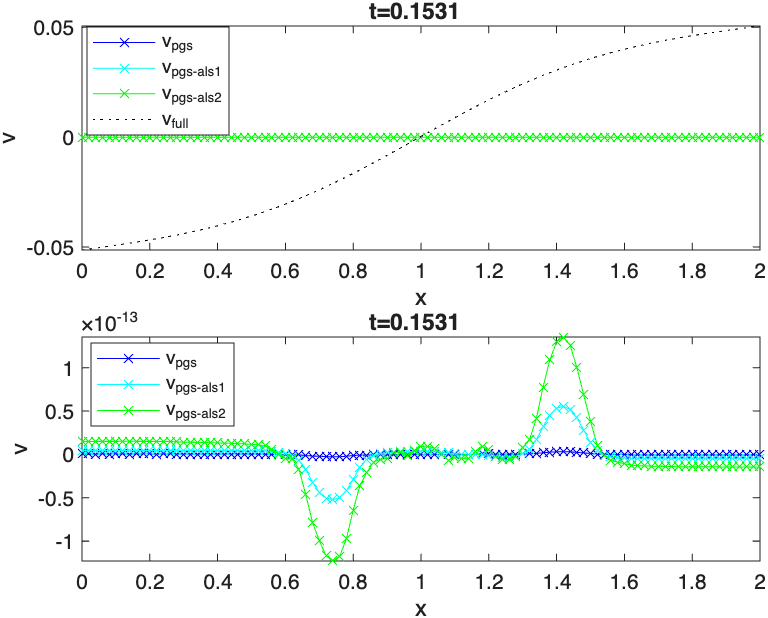}
        \end{center}
      \caption{State with Feedback Control}\label{figure2}
      }}
   \end{figure}

In Figure \ref{figure1}, we compare the states with feedback controls (\ref{linear_system_1d_a}) at various times using where the controls $u_c$ from PGS, PGS+ALS1 and  PG-Sparse+ALS2 and the full model. Figure \ref{figure2} is the zoomed plot of the bottom-right graph in Figure \ref{figure1}. It shows that at state with feedback control $u_c$ of PGS converges the fastes to $0$ at $t=0.1531$ while other state dynamics have not yet converged. The stopping criteria for these experiments is set at $\Vert v \Vert < 10^{-14}$.


 The computational complexity of the state-dependent Riccati equation is $\mathcal{O}(N_t N_x^3)$
where $N_x$ is the state dimension and $N_t$ is the number of time steps. For sparse tensor data, the approximated tensor rank $R$ is much smaller than $N_x$;i.e. $R \ll N_x$. Consequently, the time $n_t$ it takes to stabilize the reduced dynamics is much smaller than $N_t$. In Table \ref{table1}, the first column is the relative cpu times, $\frac{J_{full}-J_{red}}{J_{full}}$. The cost associated with PGS is the least among the costs. However there is an off-line cost for finding the rank $R$. The computational cost of a tensor CP decomposition using ALS or PGS is $\mathcal{O}(3IJKR)$ where $R \ll max\{I,J,K\}$.
 In conclusion, our numerical experiments show that the state dynamics with feedback control  from the PGS method is fastest with the least cost.

\begin{table}[h]
\caption{Computational Cost }
\label{table1}
\begin{center}
\begin{tabular}{|c||c|c|}
\hline
  & relative CPU time & Computational Complexity \\
\hline
$J_{full}$ & $1$ & $\mathcal{O}(N_t N_x^3)$\\
\hline
$J_{pgs}$ & $1.965 \times 10^{-29}$ & $\mathcal{O}(n_t R^3)$\\
\hline
$J_{pgs-als1}$ & $6.3809 \times 10^{-27}$  & $\mathcal{O}(n_t R^3)$\\
\hline
$J_{pgs-als12}$ & $3.896 \times 10^{-26}$ & $\mathcal{O}(n_t R^3)$\\
\hline
\end{tabular}
\end{center}
\end{table}

\section{CONCLUSIONS}

We proposed a method of approximating tensor rank from a tensor datasets of snapshots of the solution profile of a state dynamics from sampled parameters. Through CP decomposition powered by proximal gradient method for sparsity and Golub-Kahan bidiagonalization as a regularization, we find the best basis sets with an estimated low rank. The low rank estimation plays an important role since it directly dictates the reduction in the state dimension of the original full order model. The main advantages of our proposed method:  fast control feedback stabilization, huge computational savings and easily constructed tensor datasets.









\end{document}